\documentclass[11pt]{article}
\usepackage{amsfonts}
\usepackage{amssymb}
\usepackage{amssymb,amsfonts,amsmath,amsthm,cite, color}
\usepackage{epsfig}
\newtheorem{thm}{Theorem}[section]

\newtheorem{lem}{Lemma}[section]

\makeatletter \@addtoreset{equation}{section}
\def\mid{{\,|\,}}
\def\c3{{c\phi_3}}
\def\cc3{{\overline{c}\phi_3}}

\def\qed{\hfill \rule{4pt}{7pt}}

\def\Q{{\mathbb Q}}

\def\Z{\mathbb Z}

\begin{document}

%==================================================================
\begin{center}
{{\Large\bf Congruences for an arithmetic function from 3-colored Frobenius partitions }}

\vskip 6mm

{ Laizhong Song and Xinhua Xiong
\\[%
2mm] Department of Mathematics\\
China Three Gorges  University, YiChang, Hubei  443002,
P.R. China \\[3mm]
Slz@ctgu.edu.cn and XinhuaXiong@fudan.edu.cn \\[0pt%
] }
\end{center}

%=================================================================
\begin{abstract}
Let $a(n)$ defined by
\begin{displaymath}
\sum_{n=1}^{\infty}a(n)q^n := \prod_{n=1}^{\infty}\frac{1}{(1-q^{3n})(1-q^n)^3}.
\end{displaymath}
In this note, we prove that for every nonnegative integer $n$,
\begin{eqnarray*}
\begin{split}
a(15n+6)&\equiv 0\pmod{5},\\
a(15n+12)&\equiv 0\pmod{5}.\\
\end{split}
\end{eqnarray*}
As a corollary, we obtained some results of Ono\cite{Ono96}.
\end{abstract}
%=================================================================



\section{Introduction}
Arithmetic functions defined by infinite products have many beautiful congruence properties.
A typical example is  Euler's partition function $p(n)$, which is given by
\begin{displaymath}
\sum_{n=0}^{\infty}p(n)q^n=\prod_{n=1}^{\infty}\frac{1}{(1-q^n)}
\end{displaymath}
In this note, Let $a(n)$ defined by
\begin{displaymath}
\sum_{n=1}^{\infty}a(n)q^n := \prod_{n=1}^{\infty}\frac{1}{(1-q^{3n})(1-q^n)^3}.
\end{displaymath}
we prove that
\begin{thm}\label{thm1.1}
for every nonnegative integer $n$,
\begin{eqnarray*}
\begin{split}
a(15n+6)&\equiv 0\pmod{5},\\
a(15n+12)&\equiv 0\pmod{5}.\\
\end{split}
\end{eqnarray*}
\end{thm}
As a corollary of our results, we obtain the following result which is  due to Ono\cite{Ono96},
\begin{thm}\label{thm1.2}
\begin{eqnarray*}
\begin{split}
\c3(45n+7)&\equiv 0\pmod{5},\\
\c3(45n+22)&\equiv 0\pmod{5},\\
\c3(45n+37)&\equiv 0\pmod{5}.\\
\end{split}
\end{eqnarray*}
Where $\c3(n)$ is the number of three-colored Frobenius partitions of $n$.
\end{thm}

\section{Preliminary Facts}
This section contains the basic definitions and some properties of modular forms that we will
use later. For more details see Koblitz \cite{Koblitz84}.

Define $\Gamma=Sl_2(\Z),\,\Gamma_0(N):=\left\{\left(
           \begin{array}{cc}
             a & b \\
             c & d \\
           \end{array}
         \right)\mid ad-bc=1,c\equiv 0\,\,\mbox{mod}\,\,N  \right\}$. Let $M_k(\Gamma_0(N))$ denote
the vector space of holomorphic modular forms on $\Gamma_0(N)$ of weight $k$.

If $f(z)$ in $M_k(\Gamma_0(N))$ has the Fourier
expansion $f(z)= \sum_{n=0}^{\infty}a(n)q^n,$ for an integer $t\mid N$, define the U-operator by
\begin{displaymath}
f(z)\mid U(t):=\sum_{n=0}^{\infty}a(tn)q^n
\end{displaymath}
It is a standard fact that $f(z)\mid U(t) \in M_k(\Gamma_0(N)),$ if $f(z)\in M_k(\Gamma_0(N)).$ Before
starting the proof, we need some lemmas. The Dedekind Eta-function is defined by the infinite product
\begin{displaymath}
\eta(z):=q^{\frac{1}{24}}\prod_{n=1}^{\infty}(1-q^n),
\end{displaymath}
where $q=e^{2\pi iz}$ and ${\rm Im}(z)>0$. A function $f(z)$ is called an Eta-product if it can be written in the form of
\begin{displaymath}
f(z)=\prod_{\delta\mid N}\eta^{r_\delta}(\delta z)
\end{displaymath}
where $N$ and $\delta$ are positive integers and $r_{\delta}\in \mathbb{Z}$.

\begin{lem}[\cite{Gordon90}]\label{lem2.1} If $f(z)=\prod_{\delta|N}\eta^{r_{\delta}}(\delta
z)$ is an Eta-product satisfying
the following conditions:
\begin{displaymath}
  k=\frac{1}{2}\sum_{\delta|N}r_{\delta}\in 2\mathbb{Z},\,\,
\sum_{\delta|N}\delta r_{\delta}\equiv 0 \ ({\rm mod}\ 24),\,\,
\sum_{\delta|N}\frac{N}{\delta} r_{\delta}\equiv 0 \ ({\rm mod}\
24),
\end{displaymath} and  $\prod_{\delta\mid N}\delta^{r_\delta}\in \Q^2.$
then $f(z)$ satisfies
\begin{equation*}\label{relation1}
f\left(\frac{az+b}{cz+d}\right)=(cz+d)^kf(z)
\end{equation*}
for each $\left(
           \begin{array}{cc}
             a & b \\
             c & d \\
           \end{array}
         \right)\in \Gamma_0(N)$.
\end{lem}
\begin{lem}[\cite{Ligozat75}]\label{lem2.2}
Let $c,d$ and $N$ be positive integers with $d|N$ and $(c,d)=1$. If
$f(z)$ is an Eta-product satisfying the conditions in Lemma
\ref{lem2.1} for $N$, then the order of vanishing of $f(z)$ at the
cusp $\frac{c}{d}$ is
\[
\frac{N}{24}\sum_{\delta
|N}\frac{(d,\delta)^2r_{\delta}}{(d,\frac{N}{d})d\delta}.
\]
\end{lem}
\begin{lem}[\cite{Ono2}]\label{lem2.3}
Let $ f(z) = \sum_{n=0}^{\infty}a(n)q^n$ be a modular form in $M_k(\Gamma_0(N))$ and let
$d:=\,\mbox{gcd}\,(r,t).$ If $0\le r < t,$ then
\begin{displaymath}
f_{r,t}(z)=\sum_{n\equiv r\pmod{t}}a(n)q^n
\end{displaymath}
is the Fourier expansion of a modular form in $M_k(\Gamma_1(\frac{Nt^2}{d}))$
\end{lem}
In order to state the last lemma, we need

\noindent Definition \,Let $l$ be a positive integer and $f = \sum_{n \geq N} a(n) q^n$ a
formal power series in the variable $q$ with rational integer coefficients. The $l$-adic
order of $f$ is defined by
\begin{displaymath}
Ord_l(f) = \mbox{inf}\{n \mid a(n)\not\equiv 0\,\,\mbox{mod}\,\,l\}
\end{displaymath}
Sturm proved the following criterion for determining whether two modular forms are
congruent for primes , One\cite{Ono2} noted the criterion holds for general integers as modulus.
\begin{lem}[\cite{Sturm}]\label{lem2.4}
Let $l$ be a positive integer and  $f(z),\, g(z)\in M_k(\Gamma_0(N))$ with rational integers
satisfying
\begin{displaymath}
ord_l(f(z)-g(z))\geq 1+ \frac{kN}{12}\prod_{p}(1+\frac{1}{p}),
\end{displaymath}
 where the product is over
the prime divisors $p$ of $N$. Then $f(z)\equiv g(z)\ ({\rm mod\
}l)$, i.e., $ord_l(f(z)-g(z))=\infty$.
\end{lem}

\section{Proof of Theorem \ref{thm1.1}}
We construct modular forms to prove our theorem.

\noindent {\it Proof of Theorem \ref{thm1.1}.}  Define the Eta-product
\begin{eqnarray*}
f(z)&=&\frac{\eta^{7}(15z)\eta^{9}(45z)}{\eta^3(z)\eta(3z)}\nonumber\\[5pt]
&=& q^{21}\left(\prod_{n=1}^{\infty}\frac{1}{(1-q^{n})^3(1-q^{3n})}\right)
\left(\prod_{n=1}^{\infty}(1-q^{15n})^{7}(1-q^{45n})^{9}\right)\nonumber\\[5pt]
&=:& \sum_{n\geq 21} r(n)q^n
\end{eqnarray*}
It  easily can be seen that our Theorem \ref{thm1.1} is equivalent to the  two congruences: for every nonnegative integer $n$,
\begin{displaymath}
r(15n+3)\equiv 0\pmod{5} \,\,\mbox{and}\,\, r(15n+12)\equiv 0\pmod{5}.
\end{displaymath}
Setting $N=45$, by Lemma \ref{lem2.1} and \ref{lem2.2} it is easy to check that $f(z)\in M_6(\Gamma_0(45))$.   In order to prove the equivalent congruences above,
We use Lemma \ref{lem2.3}and \ref{lem2.4} To be more precise, define a modular form
\begin{displaymath}
f_{3, 15}(z):= \sum_{n\equiv 3\pmod{15}, n\geq 21}r(n)q^n
\end{displaymath}
By Lemma \ref{lem2.3} we find that $f_{3, 15}(z)\in M_6(\Gamma_1(3375)).$ By Lemma \ref{lem2.4}, $r(15n+3)\equiv 0\pmod{5}$ for every
nonnegative integer $n$ can be proved if it holds for the first 12150001 terms which is  verified by machinery calculation. The another congruence
is proved similarly.\qed

\section{Proof of theorem \ref{thm1.2}}
Recall that the three-colored partition function $\c3(n)$ has generating function
\begin{displaymath}
\sum_{n=0}^{\infty}\c3(n)q^n=\sum_{n=0}^{\infty}p(n)^n + 9q \prod_{n=1}^{\infty} \frac{(1-q^{9n})^3}{(1-q^{3n})(1-q^n)^3}:=
\sum_{n=0}^{\infty}(p(n)+\cc3(n))q^n
\end{displaymath}
Where $p(n)$ is the ordinary partition function with the convention that $p(n)=0$, if $n\notin \Z.$ We will use modified Jacobi formula
\begin{displaymath}
q\prod_{n=1}^{\infty}(1-q^{9n})^3=\sum_{n=0}^{\infty}(-1)^n(2n+1)q^{(9n^2+9n)/2+1}.
\end{displaymath}
\noindent {\it Proof of Theorem \ref{thm1.2}.} By Jacobi formula above,
\begin{eqnarray*}
\c3(45n+7)&=& \cc3(45n+7)\\
           &= &9\sum_{k\geq 0}(-1)^k(2k+1)(a(45n+7-(1+\frac{9k^2+9k}{2}))),\\
\c3(45n+22)&=& \cc3(45n+22)\\
           &= &9\sum_{k\geq 0}(-1)^k(2k+1)(a(45n+22-(1+\frac{9k^2+9k}{2}))),\\
\c3(45n+37)&=& \cc3(45n+22)\\
           &= &9\sum_{k\geq 0}(-1)^k(2k+1)(a(45n+37-(1+\frac{9k^2+9k}{2}))),\\
\end{eqnarray*}
Modulo $45$, $1+(9k^2+9k)/2$ is $1,\, 10$, or $28$ and therefore $45n+7-(1+\frac{9k^2+9k}{2}),\, 45n+22-(1+\frac{9k^2+9k}{2})$ and $45n+37-(1+\frac{9k^2+9k}{2})$ are congruent to $45n+6,\, 9, \,12, \,21, \,24, \,27, \,36,\, 39, \,42.$  The cases of $45n+9, 24, 39$ correspond to $ 1+(9k^2+9k)/2\equiv 0\pmod{45}$, which is equivalent to $k\equiv 2\pmod{5}$, we find that in this case the factors $2k+1\equiv 0\pmod{5}$, Moreover $45n+6,\,  12,\, 21, \, 27,\, 36,\, 42$ are always the form of $15n+6$ or $15n+12$. So by Theorem \ref{thm1.1}, Theorem \ref{thm1.2} is proved. \qed

%-----------------------------------------------------------------

\end{document}